\documentclass[a4paper,11pt]{amsart}

\usepackage{graphicx}
\usepackage{amsmath}
\usepackage{amsthm}
\usepackage{amssymb}
\usepackage{mathtools}
\usepackage{tikz-cd}
\usepackage{enumitem}
\usepackage{hyperref}
\usepackage{url}

\numberwithin{equation}{section}
\setenumerate{fullwidth,label={\rm(\arabic*)}}

\theoremstyle{plain}
\newtheorem{mainthm}{Theorem}

\newtheorem{thm}{Theorem}[section]
\newtheorem{prop}[thm]{Proposition}

\newtheorem{lem}[thm]{Lemma}
\theoremstyle{definition}

\newtheorem{exam}[thm]{Example}
\theoremstyle{remark}

\DeclareMathOperator{\Hom}{Hom}
\DeclareMathOperator{\uHom}{\underline{Hom}}

\DeclareMathOperator{\Ext}{Ext}

\DeclareMathOperator{\pd}{pd}
\DeclareMathOperator{\id}{id}
\DeclareMathOperator{\gd}{gl.dim}

\DeclareMathOperator{\rad}{rad}

\DeclareMathOperator{\add}{add}

\renewcommand{\mod}{\operatorname{mod}}
\DeclareMathOperator{\umod}{\underline{mod}}

\DeclareMathOperator{\proj}{proj}

\newcommand{\la}{\langle}
\newcommand{\ra}{\rangle}
\mathchardef\mhyphen="2D

\newcommand{\bo}{\mathrm{b}}
\newcommand{\mb}{\mathrm{-,b}}

\newcommand{\sg}{\mathrm{sg}}

\newcommand{\bfD}{\mathbf{D}}

\newcommand{\bfK}{\mathbf{K}}

\title{A note on singular equivalences and idempotents}
\author{Dawei Shen}
\address{School of Mathematics and Statistics \\
         Henan University \\
         475004, Kaifeng, Henan \\
         P. R. China}
\email{sdw12345@mail.ustc.edu.cn}

\begin{document}
\begin{abstract}
Let $\Lambda$ be an Artin algebra and let $e$ be an idempotent in $\Lambda$.
We study certain  functors which preserve the singularity categories.
Suppose $\pd \Lambda e_{e\Lambda e}<\infty$  and
$\id{}_\Lambda\tfrac{\Lambda/\la e\ra}{\rad \Lambda/\la e\ra} < \infty$,
we show that there is a singular equivalence between $e\Lambda e$ and $\Lambda$.
\end{abstract}

\maketitle

\section{Introduction}
Let $R$ be a left noetherian ring.
Let $\mod R$ be the category of finitely generated left $R$-modules and
let $\proj R$ be the full subcategory of projective left $R$-modules.
A complex  in $\mod R$ is called {\em perfect} if it is quasi-isomorphic to a bounded complex in $\proj R$.
The {\em singularity category} of $R$ is the Verdier quotient of the bounded derived category
of $\mod R$ by the thick subcategory of perfect complexes  ~\cite{Buc1987, Orl2004}.
If there is a triangle equivalence between the singularity categories of two  rings $R$ and $S$,
then such an equivalence is called a {\em singular equivalence} ~\cite{Chen2011}.

Let  $e$ be an idempotent  in $R$. There is a recollement
\[\begin{tikzcd}
\mod R/\la e\ra  \arrow{rr} &&
\mod R \ar[shift left=2.5]{ll} \ar[shift right=3]{ll} \ar{rr}[description]{i} &&
\mod eRe \ar[shift left=2.5]{ll} \ar[shift right=3]{ll}[description]{i_\lambda}.
\end{tikzcd}\]
Here, $i$ takes any  $M$ in $\mod R$ to $eM$ and $i_\lambda$ takes any $N$ in $\mod eRe$ to $Re\otimes_{eRe}N$.
Suppose
\begin{enumerate}
   \item $\pd{}_{eRe}eR<\infty$, and
   \item  $\pd{}_RM<\infty$ for every  $M\in \mod R$ annihilated by $e$,
\end{enumerate}
then $i$ induces a singular equivalence between $R$ and $eRe$; see \cite{Chen2009}.

In the present paper, we investigate the singular equivalence induced by the functor $i_\lambda$.
We have the following

\begin{mainthm} \label{thm:A}
Let $R$ be a left noetherian ring and $e$ be an idempotent in $R$. Suppose
\begin{enumerate}
  \item $\pd Re_{eRe}<\infty$, and
  \item every $M\in\mod R$ admits a projective resolution $P$ such that $P^{-i}\in \add Re$ for  every sufficiently large $i$,
\end{enumerate}
then $i_\lambda$ induces singular equivalence between $eRe$ and $R$.
\end{mainthm}

Let $\Lambda$ be an Artin algebra and $e$ be an idempotent in $\Lambda$.
The following conditions  are studied in ~\cite{PSS2014}.
\[\begin{split}
(\alpha)\,  \id{}_\Lambda\tfrac{\Lambda/\la e\ra}{\rad \Lambda/\la e\ra} < \infty & \qquad (\beta)  \,  \pd{}_{e\Lambda e} e\Lambda  < \infty \\
(\gamma)\,  \pd{}_\Lambda\tfrac{\Lambda/\la e\ra}{\rad \Lambda/\la e\ra} < \infty & \qquad  (\delta)\,  \pd \Lambda e_{e\Lambda e} < \infty
\end{split}\]
It turns out that
\begin{enumerate}
  \item $(\alpha)$ and  $(\beta)$ hold if and only if  $(\gamma)$ and  $(\delta)$ hold;
  \item $i$ induces a singular equivalence  if and only if  $(\beta)$ and $(\gamma)$ hold.
\end{enumerate}

As a complement, we have the following

\begin{mainthm} \label{thm:B}
Let $\Lambda$ be an Artin algebra and $e$ be an idempotent in $\Lambda$.
Then $i_\lambda$ induces a singular equivalence if and only if  $(\alpha)$ and  $(\delta)$ hold.
\end{mainthm}

\section{Singularity categories}
Let $R$ be a left noetherian ring.
Let $\bfK^\mb(\proj R)$ be the  homotopy category of
bounded above complexes in $\proj R$ with bounded cohomologies.
It is a triangulated category, whose translation functor $\Sigma$ is the shift of complexes ~\cite{Ver1977}.

A complex $P$ in $\bfK^\mb(\proj R)$ is  perfect if there is an integer $\ell$ such that
the $(-i)$-th coboundary of $P$ is projective for every  $i\geq\ell$.
Let $\bfK^\bo(\proj R)$ be the full subcategory of perfect complexes.

The  singularity category of $R$ is the quotient triangulated category
\[\bfD_\sg(R)\coloneqq  \bfK^\mb(\proj R)/\bfK^\bo(\proj R).\]

Let $\umod R$ be the projectively stable category of $\mod R$.
For any $M$ in $\mod R$, take a projective precover $\pi\colon P\to M$;
the syzygy $\Omega M$ of $M$ is the kernel of $\pi$.
There is a functor
\[\Omega \colon \umod R\to \umod R\]
sending a module $M$ to the syzygy $\Omega M$.
For any $M$ in  $\mod R$, take a projective resolution $pM$ of $M$.
There is a functor
\[p\colon \umod R\to \bfD_\sg(R)\]
sending a module $M$ to the projective resolution $pM$.
For any  $\ell \geq 0$, there is a natural isomorphism
\[ p\Omega^\ell (M)\cong\Sigma^{-\ell} pM.\]

\begin{lem}[{\cite[Lemma~2.1]{Chen2011}}] \label{lem:singcat1}
For any $X \in \bfD_\sg(R)$,
there is  $n\geq 0$ and  $M\in\mod R$ such that
\[\Sigma^{-n}X\cong pM.\]
\end{lem}

\begin{lem}[{\cite[Exemple~2.3]{KV1987}}]   \label{lem:singcat2}
For any $M,N\in \mod R$, there is an isomorphism
\[\underrightarrow{\lim}_{\ell\geq 0}\uHom_R(\Omega^\ell M,\Omega^\ell N)\cong \Hom_{\bfD_\sg(R)}(pM,pN).\]
\end{lem}

Let  $e$ be an idempotent in $R$.
There is a functor
\[i_\lambda \colon \mod eRe \to \mod R\]
such that $i_\lambda(M)=Re\otimes_{eRe}M$ for every $M\in\mod eRe$.
It is a full faithful  functor which takes projectives to projectives.
The functor $i_\lambda$ restricts to a full faithful functor $i_\lambda'$ such that
\[i_\lambda'\colon \proj eRe \overset{\simeq}\to \add Re  \overset{\subset}\to  \proj R.\]
Here, we denote by  $\add Re$  the full subcategory  of summands
of finite direct  sums of copies of $Re$.

We now restate and  prove Theorem ~\ref{thm:A}.

\begin{thm}
Suppose  $\pd Re_{eRe}<\infty$, then $i_\lambda$ induces a full faithful triangle functor
\[ \bfD_\sg(i_\lambda)\colon \bfD_\sg(eRe) \to \bfD_\sg(R).\]
Moreover, $\bfD_\sg(i_\lambda)$ is a triangle equivalence if and only if for every $M\in\mod R$,
there is an integer $\ell$ and a projective resolution $P$ of $M$  such that $P^{-i}\in \add Re$ for every $i\geq \ell$.
\end{thm}

\begin{proof}
Let  $X$ be  in $\bfK^\mb(\proj eRe)$ which is exact at  degree $\leq \ell$.
Then $i_\lambda(X)$ is exact at degree  $\leq \ell-\pd Re_{eRe}$.
Therefore $i_\lambda$ induces a triangle functor
\[\bfK^\mb(i_\lambda)\colon \bfK^\mb(\proj eRe) \to \bfK^\mb(\proj R).\]
Since $i_\lambda$ is full faithful, $\bfK^\mb(i_\lambda)$ is also  full faithful.
Since $i_\lambda$ preserves perfect complexes,
it induces a triangle functor $\bfD_\sg(i_\lambda)$ such that
\[\bfD_\sg(i_\lambda)\colon \bfD_\sg(eRe) \overset{\simeq}\to \bfK^\mb(\add Re)/\bfK^\bo(\add Re) \to \bfD_\sg(R).\]

Let $f\colon X\to Y$ be a morphism in $\bfK^\mb(\proj R)$,
where $X^i$ is zero for every $i\leq \ell$ and $Y$ belongs to $\bfK^\mb(\add Re)$.
Let $\sigma_{\geq \ell}Y$ be the stupid truncation of $Y$ at degree $\geq \ell$,
then it belongs to  $\bfK^\bo(\add Re)$.
Since $f$ factors through  $\sigma_{\geq \ell}Y$,
by ~\cite[Proposition 10.2.6]{KS2005} the functor $\bfD_\sg(i_\lambda)$ is full faithful.

By Lemma ~\ref{lem:singcat1} and ~\ref{lem:singcat2},
the singularity category of $R$ is triangle equivalent to the stabilization of the stable category $\umod R$.
Then the denseness of  $\bfD_\sg(i_\lambda)$ follows from ~\cite[Corollary ~2.13]{Chen2018}.
\end{proof}

\section{Triangular matrix rings}
Let $T$ and $S$ be two rings and  $M$ be an $S$-$T$-bimodule.
We consider the triangular matrix ring
\[R=\begin{pmatrix}
      T & 0 \\
      M & S
\end{pmatrix}.\]
Following ~\cite[III.2]{ARS1995} a left $R$-module  is  given by
\[(X,Y,\phi)=\left\{\begin{pmatrix}x \\
y\end{pmatrix}\mid x\in X, y\in Y\right\},\]
where $X$ is a left $T$-module, $Y$ is a left $S$-module and
$\phi\colon M\otimes_TX\to Y$ is a left $S$-module map.
The action is given by
\[\begin{pmatrix}
      t & 0 \\
      m & s
\end{pmatrix}\begin{pmatrix}
      x \\
      y
\end{pmatrix} =\begin{pmatrix}
      tx \\
      \phi(m\otimes x)+sy
\end{pmatrix}\]
for every $t\in T$, $s\in S$ and $m\in M$.

Let $e=\mathrm{diag}(0,1)$ be  in $R$.
Then  $Re=eRe$ and
\[\add Re=\{(0,Y,0)\mid Y\in \proj S\}.\]

\begin{lem} \label{lem:tri}
Let $X$ be a left $T$-module and $Y$ be a left $S$-module.
\begin{enumerate}
\item $\pd{}_SY=\pd{}_R(0,Y,0)$;
\item $(X,M\otimes_T X,1)$ is a projective left  $R$-module if and only if $X$ is a projective left $T$-module.
\end{enumerate}
\end{lem}

We have the following; compare ~\cite[Theorem 4.1]{Chen2009}.

\begin{prop}
Let $T$ and $S$ be left noetherian rings,
and let $M$ be an $S$-$T$-bimodule such that ${}_SM$ is finitely generated.
Assume that $T$   has finite left global dimension, then there is a triangle equivalence
\[\bfD_\sg(\begin{pmatrix}
      T & 0 \\
      M & S\end{pmatrix})
\cong \bfD_\sg(S). \]
\end{prop}

\begin{proof}
One checks that the triangular matrix ring $R$ is  left noetherian.
Let $(X,Y,\phi)$ be in $\mod R$.
Since $\pd{}_TX$ is finite,
by Lemma ~\ref{lem:tri} there a projective resolution $P$ of $(X,Y,\phi)$ such that
$P^{-i}\in \add Re$ for every $i>\pd{}_TX$.
Then by Theorem \ref{thm:A} there is a triangle equivalence between $\bfD_\sg(R)$ and $\bfD_\sg(S)$.
\end{proof}

\begin{exam}[see  ~\cite{Chen2011}, Proposition 4.1]
Let $k$ be a field, $\Lambda$ be a finite dimensional $k$-algebra and $M$ be a finite dimensional left $\Lambda$-module.
Then $M$ is a $\Lambda$-$k$-bimodule.
The {\em one-point extension} of $\Lambda$ by $M$ is the  triangular matrix algebra
\[\Lambda[M]=\begin{pmatrix}
      k & 0 \\
      M & \Lambda\end{pmatrix}.\]
By Theorem ~\ref{thm:A} there a singular equivalence between $\Lambda$ and $\Lambda[M]$.
\end{exam}

\section{Artin algebras}
Let $\Lambda$ be an Artin algebra. We need the following well known

\begin{lem} \label{lem:ext-hom}
Let $M$ be in $\mod \Lambda$ and $P_M$ be a minimal projective resolution of $M$.
For any semi-simple $\Lambda$-module  $S$ and any  $i\geq 0$ there is an isomorphism
\[\Ext^i_\Lambda(M,S)\cong \Hom_\Lambda(P_M^{-i},S).\]
\end{lem}

\begin{proof}
For $i= 0$, it is obvious.
For $i = 1$, let  $K$ be the $(-1)$-th coboundary of $P_M$.
There is an exact sequence
\[ K \rightarrowtail P_M^0 \overset{\pi} \twoheadrightarrow  M.\]
Applying $\Hom_\Lambda(-,S)$ to it,
we obtain an exact sequence
\[ \Hom_\Lambda(M,S) \overset{\pi^*} \rightarrowtail \Hom_\Lambda(P_M^0,S)\to \Hom_\Lambda(K,S)\twoheadrightarrow \Ext^1_\Lambda(M,S). \]
Since $\pi$ is a projective cover and $S$ is semi-simple, ${\pi^*}$ is surjective.
Then  \[\Hom_\Lambda(K,S) \cong \Ext^1_\Lambda(M,S).\]
Since $P_M^{-1}$ is a projective cover of $K$, there is an isomorphism
\[\Ext^1_\Lambda(M,S)\cong \Hom_\Lambda(P_M^{-1},S).\]
By shifting one proves the isomorphism for $i\geq 2$.
\end{proof}

For any $P\in \proj \Lambda$, we have the following
\begin{equation}\label{eq:eq1}
P \in \add \Lambda e \iff \Hom_\Lambda(P,\tfrac{\Lambda/\la e\ra}{\rad \Lambda/\la e\ra})=0.
\end{equation}

We now restate and  prove Theorem ~\ref{thm:B}.

\begin{thm}
Let $\Lambda$ be an Artin algebra and $e$ be an idempotent in $\Lambda$.
Suppose $\pd \Lambda e_{e\Lambda e}<\infty$, then $\bfD_\sg(i_\lambda)$ is a triangle equivalence
if and only if $\id{}_\Lambda\tfrac{\Lambda/\la e\ra}{\rad \Lambda/\la e\ra}<\infty$.
\end{thm}

\begin{proof}
$``\Longrightarrow"$
Let $P$ be a minimal projective resolution of $\Lambda/\rad \Lambda$.
If $\bfD_\sg(i_\lambda)$ is dense, by Theorem ~\ref{thm:A}
there is  $\ell\geq 0$ such that $P^{-i} \in \add \Lambda e$ for every $i>\ell$.
By Lemma \ref{lem:ext-hom} and \eqref{eq:eq1} we have
\[\Ext^i_\Lambda(\Lambda/\rad \Lambda,\tfrac{\Lambda/\la e\ra}{\rad \Lambda/\la e\ra})=0.\]
Then $\id{}_\Lambda\tfrac{\Lambda/\la e\ra}{\rad \Lambda/\la e\ra}$ is finite.

$``\Longleftarrow"$
Let $M$ be  in  $\mod \Lambda$
and  $P_M$ be a minimal projective resolution  of $M$.
If $\id{}_\Lambda\tfrac{\Lambda/\la e\ra}{\rad \Lambda/\la e\ra}$ is finite,
by  Lemma ~\ref{lem:ext-hom} we have
\[P_M^{-i}\in\add \Lambda e,\forall\, i>\id{}_\Lambda\tfrac{\Lambda/\la e\ra}{\rad \Lambda/\la e\ra}.\]
From Theorem ~\ref{thm:A} we infer that $\bfD_\sg(i_\lambda)$ is a triangle equivalence.
\end{proof}

We end this section by an example.

Let $\Lambda$ be an Artin algebra.
Let $S$ be a semi-simple left $\Lambda$-module with $\id_\Lambda S \leq 1$.
Denote by
\[{}^\perp S = \{M\in\mod \Lambda\mid \Hom_\Lambda(M,S)=\Ext^1_\Lambda(M,S)=0\}\]
the  perpendicular category of $S$ in $\mod \Lambda$.
By Lemma \ref{lem:ext-hom} a finitely generated left $\Lambda$-module $M$ belongs to ${}^\perp S$  if and only if
$M$ admits a projective presentation
\[P^{-1}\to P^0\twoheadrightarrow M\]
such that both $P^{-1}$ and $P^0$ belong to  ${}^\perp S$.

Let $e$ be an idempotent in $\Lambda$ such that  $\proj\Lambda\cap {}^\perp S$ coincides with $\add \Lambda e$.
Recall the functor $i_\lambda$ decomposes into
\[i_\lambda \colon \mod e\Lambda e \overset{\simeq}\to {}^\perp S \overset{\subseteq}\to \mod \Lambda.\]

We have the following; compare ~\cite[Proposition ~2.13]{CY2013}.

\begin{prop}
Keep the notation as previous.
\begin{enumerate}
  \item $\gd e \Lambda e\leq \gd \Lambda \leq \gd e\Lambda e+2$;
  \item there is a singular equivalence between $e\Lambda e$ and $\Lambda$.
\end{enumerate}
\end{prop}

\begin{proof}
(1)  Let $\Omega^2M$ be the minimal second syzygy for $M$ in $\mod \Lambda$.
Since $\id_\Lambda S\leq 1$, $\Ext_\Lambda^i(M,S)=0$ for every $i\geq 2$.
By Lemma ~\ref{lem:ext-hom}  $\Omega^2M \in {}^\perp S$.
Since $i_\lambda$ preserves projective resolutions, for any $N$ in $\mod e\Lambda e$ we have
\[\pd_{e \Lambda e} N= \pd_\Lambda i_\lambda(N).\]
Then
\[\pd_\Lambda M \leq \pd_\Lambda\Omega^2M+2 = \pd_{e \Lambda e} i(\Omega^2M)+2.\]

(2) From ~\cite[Proposition ~1.1]{GL1991} we infer that  ${}^\perp S$ is an exact subcategory of $\mod \Lambda$.
Then the right ${e\Lambda e}$-module $\Lambda e$ is  projective.
By Theorem ~\ref{thm:B} there is a singular equivalence between $e\Lambda e$ and $\Lambda$.
\end{proof}

\section*{Acknowledgments}
The author is thankful to Professor Xiao-Wu Chen for discussions.
This work is supported by  National Natural Science Foundation of China  (No. 11801141).

\end{document}